\documentstyle[12pt]{article}
\begin{document}

\title{ \Large {The Geometry and Topology on  Grassmann Manifolds}
 \footnotetext{
\mbox{}\quad  \ \ {\it Subject classification}. Primary 14M15;
Secondary 53C42, 57R70. }}
\author{  Zhou Jianwei}
\date{\small Department of Mathematics, Suzhou University, Suzhou
215006, P. R. China }

\maketitle
\begin{abstract}   This paper shows that the Grassmann Manifolds $G_{\bf F}(n,N)$
can all be imbedded in an Euclidean space $M_{\bf F}(N)$ naturally
and the imbedding can be realized by  the eigenfunctions of
Laplacian $\triangle$ on $G_{\bf F}(n,N)$. They are all minimal
submanifolds in some spheres of $M_{\bf F}(N)$ respectively. Using
these imbeddings, we construct some degenerate Morse functions on
Grassmann Manifolds, show that the homology of the complex and
quaternion Grassmann Manifolds can be computed easily.
\end{abstract}

\baselineskip 15pt
\parskip 5pt

\vskip 1cm

\centerline{\bf \S 1. Introduction}

\vskip 0.3cm

 Let $G_{\bf F}(n, N)$
be the Grassmann manifold formed by all $n$-subspaces in ${\bf
F}^N$, where ${\bf F}$ is the set of real numbers, complex numbers
or quaternions. The manifold $G_{\bf F}(n, N)$ is a symmetric
space (see [7] or [8]). The Grassmann manifolds are important in
the study of the geometry and the topology, especially in the
theory of fibre bundles.

 Let $\widetilde G(n, N)$ be the
oriented Grassmann manifold formed by all oriented $n$ dimensional
subspaces of ${\bf R}^N$. In [3], Chen showed that $\widetilde
G(n, N)$ can be imbedded in the unit sphere of wedge product space
$\bigwedge^n({\bf R}^N)$ as a minimal submanifold. Takahashi  [10]
proved that a compact homogenous Riemannian manifold with
irreducible linear isotropy group admits a minimal immersion into
an Euclidean sphere,  see also Takeuchi and Kobayashi [11].

Let $M_{\bf F}(N)$ be the set of $N\times N$ matrices $A$ with
values in ${\bf F}$ such that $\overline A^t =A$. \ $M_{\bf F}(N)$
is an Euclidean space. Let $M_{\bf F}(n,N) =\{ A \in M_{\bf F}(N)
\ | \ A^2=A, \ \mbox {r} (A)=n\}$ be a subspace of $M_{\bf F}(N)$,
where $\mbox {r} (A)$ be the rank of the matrix $A$. The matrix $A
\in M_{\bf F}(n,N)$ can be viewed as  a projection on Euclidean
space ${\bf F}^N$.

For any $\pi\in G_{\bf F}(n,N)$, let $e_1,\cdots,e_n$ be an
orthonormal basis of $\pi$. Then $(e_1,\cdots,e_n)$ is an $N\times
n$ matrix. Define
$$\varphi(\pi)=(e_1,\cdots,e_n)\overline
{(e_1,\cdots,e_n)}^t=\sum\limits_i e_i \bar e_i^t.$$ We show in \S
2, the map $\varphi: \ G_{\bf F}(n,N) \to  M_{\bf F}(N)$ is an
imbedding and we have $\varphi( G_{\bf F}(n,N))= M_{\bf F}(n,N)$.
Then
$$\bigcup\limits_{n=0}^N \varphi(G_{\bf F}(n,N)) =\{ A\in M_{\bf F}(N) \ | \
A^2=A \}.$$ Let $\{ A\in M_{\bf F}(N) \ | \ \mbox{tr} A=n\}$ be a
hyperplane in $M_{\bf F}(N)$ and $S(\sqrt n)$ the sphere of
$M_{\bf F}(N)$ with radius $\sqrt n$. In \S 2, we also show that
$M_{\bf F}(n,N)$ is a minimal submanifold in the sphere $S(\sqrt
n)\bigcap \{ A\in M_{\bf F}(N) \    | \ \mbox{tr} A=n\}$. These
minimal submanifolds are the natural generalization of the famous
Veronese surface.

Let $G_{\bf F}(N)$ be the group which preserving the inner product
on Euclidean space ${\bf F}^N$. With the spaces $M_{\bf F}(n,N)$,
we can show that the Grassmann manifold $G_{\bf F}(n,N)$ can be
imbedded in the group $G_{\bf F}(N)$.

In \S 3, we  construct some degenerate Morse functions on
Grassmann Manifolds. Show that the Poincar\'{e} polynomial of
$G_{\bf F}(n,N)$ can be represented by
$$P_t(G_{\bf F}(n,N))=P_t(G_{\bf F}(n,N-1)) +t^{c(N-n)}P_t(G_{\bf F}(n-1,N-1)),$$
where ${\bf F}={\bf C}$ or ${\bf F}={\bf H}$ and $c$  the
dimension of ${\bf F}$. Then the homology of the complex and
quaternion Grassmann Manifolds can be computed easily in low
dimensional cases.

These results are consistent with the results computed by using
Schubert variety. As in [4] or [5], we consider the case of ${\bf
F}={\bf C}$. Let
$$0\leq a_0\leq a_1\leq \cdots \leq a_{n}\leq N-n$$ be a sequence of
integers. There is a natural one-one  correspondence   between the
set of $(a_0, a_1, \cdots, a_{n})$ and the generators of the
homology $H_*(G_{\bf C}(n,N))$. The dimension of $(a_0, a_1,
\cdots, a_{n})$ is $2(a_0+ a_1+ \cdots + a_{n})$. Such elements
can be divided into two classes:

(1) $(a_0, a_1, \cdots, a_{n})$, where $a_n \leq N-n-1$;

(2) $(a_0, a_1, \cdots, a_{n-1}, N-n)$, where $a_{n-1} \leq N-n$.

\noindent These also show
$$P_t(G_{\bf C}(n,N))=P_t(G_{\bf C}(n,N-1)) +t^{2(N-n)}P_t(G_{\bf C}(n-1,N-1)).$$

 The Poincar\'{e} polynomial of
$G_{\bf F}(n,N)$ can also be represented by \begin{eqnarray*} &
P_t(G_{\bf F}(n,N))= & t^{cn}P_t(G_{\bf F}(n,N-2)) +
t^{c(N-n)}P_t(G_{\bf F}(n-2,N-2))\\
& &  + (1+ t^{c(N-1)})P_t(G_{\bf F}(n-1,N-2)),
\end{eqnarray*} where $n, N-n\geq 2, \ {\bf F}={\bf C}$ or ${\bf H}$.

\vskip 1cm \centerline{\bf \S 2. The minimal imbedding of $G_{\bf
F}(n, N)$ in the sphere}

\vskip 0.3cm

Let ${\bf F}$ be the set of real numbers ${\bf R}$, complex
numbers ${\bf C}$ or quaternions ${\bf H}$. The quaternions $ {\bf
H}$ is generated by $i,j,k=ij$. For any $u\in {{\bf F}}$,
$\overline u$ is the conjugation of $u$ ($\overline{u\cdot
v}=\overline v\cdot \overline u$ if $u, v \in {\bf H}$). For any
$\lambda\in   {\bf F}, \ \lambda$ acts on the right of
$u=(u_1,\cdots,u_N)^t\in {{\bf F}}^N$. For any
$u=(u_1,\cdots,u_N)^t, \ v=(v_1,\cdots,v_N)^t\in { {\bf F}}^N$,
$$(u,v) =\bar v^t\cdot u=\sum\limits_{A} \bar v_Au_A$$
defines an inner product on ${\bf F}^N$. Let $G_{\bf F}(N)$ be the
group acting on the left of ${\bf F}^N$ which preserving the inner
product $( , )$ on $ {\bf F}^N$. If ${\bf F}={\bf R}, \ G_{\bf
F}(N)=O(N)$ is the orthogonal group; if ${\bf F}={\bf C},\ G_{\bf
F}(N)=U(N)$ is the complex unitary group; if $ {\bf F}= {\bf H}, \
G_{\bf F}(N)=Sp(N)$ is the symplectic group.

Let  $G_{\bf F}(n, N)\approx \frac {G_{\bf F}(N)}{G_{\bf
F}(n)\times G_{\bf F}(N-n)}$   be the Grassmann manifold formed by
all subspaces in ${\bf F}^N$ of dimension $n$. Let
$e_1,\cdots,e_n, e_{n+1},\cdots,e_N$ be orthonormal frame fields
on $G_{\bf F}(n,N)$ such that the element of $G_{\bf F}(n,N)$ is
generated by $e_1,\cdots,e_n$ locally. By the method of moving
frame, there are local 1-forms $\omega^B_A$ defined by
$$de_A=\sum\limits_{B} e_B\omega^B_A, \ \omega^B_A + \overline\omega^A_B=0, \ \
A,B= 1,\cdots,N.$$ Restricting the two form $\Phi = \mbox{Re} \
(\sum\limits_{i,\alpha} \ \omega^\alpha_i\overline
\omega^\alpha_i)$ on $G_{\bf F}(n, N)$ defines a Riemannian metric
(see [4]). Unless otherwise stated, we agree on the following
arranges of the indices:
$$1\leq i,j,\cdots\leq n, \ \
n+1\leq\alpha,\beta,\cdots\leq N, \ \ 1\leq A,B,\cdots\leq N.$$

Let $M_{\bf F}(N)$ be the set of $N\times N$ matrices $A$  with
values in ${\bf F}$ such that $\overline A^t =A$. With the inner
product defined by
$$\langle A, B\rangle = \mbox{Re} \ \mbox{tr}  (AB)
=\sum\limits_{A} x_{AA} y_{AA} + 2\mbox{Re} \  \sum\limits_{A<B}
x_{AB}\bar y_{AB},$$ $A=(x_{AB}),  B=(y_{AB})\in M_{\bf F}(N)$, \
$M_{\bf F}(N)$ becomes an Euclidean space. The real dimension of
$M_{\bf F}(N)$ is $N+ \frac 12 cN(N-1)$, where $c$ is the real
dimension of ${\bf F}$.

{\bf Lemma 2.1} \ Let $e_1, \cdots,e_N$ be an orthonormal frame on
${\bf F}^N$. The following elements form an orthogonal basis of
$M_{\bf F}(N)$ with respect to the norm $\langle ,\rangle $
respectively,

(1) $e_A e_A^t, \ e_B e_C^t+ e_C e_B^t$, when ${\bf F}={\bf R}$;

(2) $e_A\bar e_A^t, \ e_B\bar e_C^t+ e_C\bar e_B^t, \ e_Bi\bar
e_C^t- e_Ci\bar e_B^t$, when ${\bf F}={\bf C}$;

(3) $e_A\bar e_A^t, \ e_B\bar e_C^t+ e_C\bar e_B^t, \ e_Bi\bar
e_C^t- e_Ci\bar e_B^t, \ e_Bj\bar e_C^t- e_Cj\bar e_B^t, \
e_Bk\bar e_C^t- e_Ck\bar e_B^t$, when ${\bf F}={\bf H}$,

\noindent where $A, B,C=1,\cdots,N, \ B<C$.

 {\bf Proof}. The proof
is a direct computation. For example, we have
$$\langle e_A\bar e_A^t,  e_B\bar e_C^t+ e_C\bar e_B^t\rangle
= \mbox{Re} \ \mbox{tr}  (e_A\bar e_C^t\delta_{AB}+ e_A\bar
e_B^t\delta_{AC})=2\delta_{AB}\delta_{AC},$$ and $$\langle e_B\bar
e_C^t+ e_C\bar e_B^t, e_Bi\bar e_C^t- e_Ci\bar e_B^t\rangle =
\mbox{Re} \ \mbox{tr} (e_Ci\bar e_C^t-e_Bi\bar e_B^t)=0. \ \ \
\Box$$

Note that the basis of $M_{\bf F}(N)$ described in Lemma 2.1 are
all real.

  For any $\pi\in G_{\bf F}(n,N)$, let $e_1,\cdots,e_n$ be an orthonormal basis of
 $\pi$. Then $(e_1,\cdots,e_n)$ is an $N\times n$ matrix.
Define
$$\varphi: \ G_{\bf F}(n,N) \to M_{\bf F}(N),$$
$$\varphi(\pi)=(e_1,\cdots,e_n)\overline
{(e_1,\cdots,e_n)}^t=\sum\limits_i e_i \bar e_i^t.$$ It is easy to
see that $\varphi(\pi)$ is independent of the choice of the
orthonormal basis $e_1,\cdots,e_n$. Let $M_{\bf F}(n,N) =\{ A \in
M_{\bf F}(N) \ | \ A^2=A, \ \mbox {r} (A)=n\}$ be a subspace of
$M_{\bf F}(N)$, where $\mbox {r} (A)$ be the rank of matrix $A$.

{\bf Lemma 2.2} \ The map $\varphi: \ G_{\bf F}(n,N) \to  M_{\bf
F}(N)$ is an imbedding and we have $\varphi( G_{\bf F}(n,N)) =
M_{\bf F}(n,N)$. The induced metric on $G_{\bf F}(n,N)$ defined by
$\varphi$ is $$2\Phi=2\mbox{Re} \ (\sum\limits_{i,\alpha} \
\omega^\alpha_i\overline \omega^\alpha_i).$$

{\bf Proof}. \ It is easy to see that $\varphi( G_{\bf
F}(n,N))\subset M_{\bf F}(n,N)$. On the other hand, the element
$A\in M_{\bf F}(n,N)$ can be viewed as a projection on ${\bf
F}^N$. Let $\pi=\{ Ax \ | \ x \in {\bf F}^N \}$ be a subspace of
${\bf F}^N$ and $e_1,\cdots,e_n$ be an orthonormal basis of $\pi$.
Therefore $Ae_i=e_i$. It is easy to see that $A=\sum\limits_i e_i
\bar e_i^t$ and
 $\varphi(\pi)=A$.
Then we can identify $M_{\bf F}(n,N)$ with $G_{\bf F}(n,N)$. Let
$e_1,\cdots,e_n, e_{n+1},\cdots,e_N$ be orthonormal frame fields
 on ${\bf F}^N$ such that $G_{\bf F}(n,N)$ is generated by $e_1,\cdots,e_n$ locally.
Hence
$$d\varphi = d\sum\limits_i e_i \bar e_i^t=
\sum\limits_{i,\alpha} e_\alpha\omega_i^\alpha \bar e_i^t +
\sum\limits_{i,\alpha} e_i\bar \omega_i^\alpha \bar e_\alpha^t.$$
We compute the case of ${\bf F}={\bf H}$ as an example, the other
cases are similar.

Let $\omega_i^\alpha = a_i^\alpha+ ib_i^\alpha +jc_i^\alpha +
kd_i^\alpha$, where $a_i^\alpha,b_i^\alpha,c_i^\alpha, d_i^\alpha$
are real 1-forms. Then \begin{eqnarray*} & d\varphi = &
\sum\limits_{i,\alpha} a_i^\alpha (e_\alpha\bar e_i^t +  e_i \bar
e_\alpha^t)  + \sum\limits_{i,\alpha} b_i^\alpha (e_\alpha i\bar
e_i^t -  e_i i\bar
e_\alpha^t) \\
& & +\sum\limits_{i,\alpha} c_i^\alpha (e_\alpha j\bar e_i^t - e_i
j \bar e_\alpha^t)  + \sum\limits_{i,\alpha} d_i^\alpha (e_\alpha
k\bar e_i^t - e_i k\bar e_\alpha^t), \end{eqnarray*} and
$$e_\alpha \bar e_i^t+e_i \bar e_\alpha^t, \
e_\alpha i \bar e_i^t-e_ii \bar e_\alpha^t, \ e_\alpha j \bar
e_i^t-e_ij \bar e_\alpha^t, \ e_\alpha k \bar e_i^t-e_ik \bar
e_\alpha^t$$ form a basis of tangent space  $TG_{\bf H}(n,N)$. By
Lemma 2.1, these vectors are orthogonal with respect to the inner
product $\langle, \rangle$. The norms of these vectors are all
$\sqrt 2$. Then
$$\langle d\varphi, d\varphi\rangle = 2\sum\limits_{i,\alpha} \ (a_i^\alpha\otimes a_i^\alpha +
b_i^\alpha\otimes b_i^\alpha +c_i^\alpha\otimes c_i^\alpha +
d_i^\alpha\otimes d_i^\alpha) =2\Phi.   \ \ \ \Box $$

For any $A\in M_{\bf F}(n,N)$, \
$$\langle A,A\rangle = \mbox{tr} \ A\overline A^t= \mbox{tr}  A= \mbox{tr}
[(e_1,\cdots,e_n)\overline {(e_1,\cdots,e_n)}^t]=
n,$$ then $|A|=\sqrt n$.
 These also show $r (A)=\mbox{tr}  A$ for any $A\in
 M_{\bf F}(N)$ with $A^2=A$. Therefore we also have
 $$M_{\bf F}(n,N) =\{ A\in M_{\bf F}(N) \ | \ A^2=A, \ \mbox{tr} A=n\},$$
and
$$\bigcup\limits_{n=0}^N \varphi(G_{\bf F}(n,N)) =\{ A\in M_{\bf F}(N) \ |
\ A^2=A \}.$$

 Let $I_N$ be the identity matrix of order
$N$. For $A\in M_{\bf F}(n,N)$, we have $(I_N-A)^2=I_N-A, \ r
(I_N-A)=\mbox{tr} (I_N-A)= N-n $, hence $I_N-A\in M_{\bf
F}(N-n,N)$. These show the map $A\to I_N-A$ gives an isometry
between the manifolds $G_{\bf F}(n,N)$ and $G_{\bf F}(N-n,N)$.

As show above, for any $A\in M_{\bf F}(n,N), \ |A|^2=\mbox{tr}
A=n$, then $M_{\bf F}(n,N)$ is in the sphere $S(\sqrt {n})=\{ B\in
M_{\bf F}(N) \ | \ |B|^2=n\}$. By $\mbox{tr}  A=n$, we know that
$M_{\bf F}(n,N)$ also in the hyperplane $\{ B\in M_{\bf F}(N) \ |
\ \mbox{tr} B=n\}$. This hyperplane can also be defined by $\{
B\in M_{\bf F}(N) \ | \ \langle B, I_N\rangle=n\}$. Then the
normal vector of this hyperplane is $I_N =\sum\limits_i e_i\bar
e_i^t + \sum\limits_\alpha e_\alpha\bar e_\alpha^t.$

When ${\bf F}={\bf R}$, any element $A\in M_{\bf R}(1,3)\approx
{\bf R}P^2$ can be represented by
$$A=\left( \begin{array}{ccc}
x_1^2 & x_1x_2 & x_1x_3 \\ x_2x_1 & x_2^2 & x_2x_3 \\ x_3x_1 &
x_3x_2 & x_3^2  \end{array} \right), \ \ x_1^2+x_2^2+x_3^2=1.$$
The map
$$B=(x_{AB})\in M_{\bf R}(3) \to
(x_{11},x_{22},x_{33}, \sqrt 2x_{12}, \sqrt 2x_{13}, \sqrt
2x_{23})\in {\bf R}^6$$ gives an isometry between these two
Euclidean spaces. Then $M_{\bf R}(1,3)$ is the famous Veronese
surface.

{\bf Theorem 2.3} \ The manifold $M_{\bf F}(n,N)$ is a minimal
submanifold in the sphere $S(\sqrt {n})\bigcap \{ B\in M_{\bf
F}(N) \ | \  \mbox{tr}  B=n\}$.

{\bf Proof}. \  With the notations used
above,
$$d\varphi =\sum e_\alpha\omega_i^\alpha \bar e_i^t + \sum e_i\bar
\omega_i^\alpha \bar e_\alpha^t,$$
$$d^2\varphi  =\cdots +
 \sum [e_j\omega_\alpha^j\omega_i^\alpha \bar e_i^t +
 e_\alpha\omega_i^\alpha\bar\omega_i^\beta \bar e_\beta^t +
 e_\beta\omega_i^\beta\bar\omega_i^\alpha \bar e_\alpha^t +
 e_i\bar\omega_i^\alpha\bar\omega_\alpha^j \bar e_j^t],
$$
where $``\cdots"$ is the part of $d^2\varphi$  which tangent to
$M_{\bf F}(n,N)$. Then the second fundamental form of the
imbedding $\varphi$ is
$$II=-\sum [e_j\bar\omega_j^\alpha\omega_i^\alpha \bar e_i^t +
 e_i\bar\omega_i^\alpha\omega_j^\alpha \bar e_j^t]+\sum [
 e_\alpha\omega_i^\alpha\bar\omega_i^\beta \bar e_\beta^t +
 e_\beta\omega_i^\beta\bar\omega_i^\alpha \bar e_\alpha^t ].
$$
The mean curvature vector  is
$$H=- \frac {N-n}{n(N-n)} \sum\limits_i  e_i \bar e_i^t +\frac n{n(N-n)} \sum\limits_\alpha
 e_\alpha\bar e_\alpha^t.$$

On the other hand, $\sum\limits_i e_i\bar e_i^t$ and $I_N =
\sum\limits_i e_i\bar e_i^t + \sum\limits_\alpha e_\alpha\bar
e_\alpha^t$ are
 the normal vectors on $S(\sqrt {n})$ and
$\{ B\in M_{\bf F}(N) \ | \  \mbox{tr}  B=n\}$ at $\sum\limits_i
e_i\bar e_i^t$
 respectively. These show $M_{\bf F}(n,N)$ is a minimal submanifold
in the sphere $S(\sqrt {n})\bigcap \{ B\in M_{\bf F}(N) \ | \
\mbox{tr} B=n\}$. \ \ \ $\Box$

The radius of  the sphere $S(\sqrt {n})\bigcap \{ B\in M_{\bf
F}(N) \ | \ \mbox{tr} B=n\}$ is $\sqrt {\frac {n(N-n)}N}$.

The above proof also shows, the second fundamental form of $M_{\bf
F}(n,N)$ has
 constant length [6]. The isometry group $G_{\bf F}(N)$  acts
on the Grassmann manifold $G_{\bf F}(n,N)$  naturally. For any
$g\in G_{\bf F}(N), \ A\in M_{\bf F}(N), \ Ad(g)A= gA\bar g^t$
defines an action of $G_{\bf F}(N)$ on $M_{\bf F}(N)$.
Furthermore, the following diagram is commutative
$$\begin{array}{ccc}
G_{\bf F}(n,N) & \stackrel {\quad \varphi \quad }
{\longrightarrow} &
M_{\bf F}(N) \\ g\downarrow  \ \ & \quad & \qquad\downarrow {Ad(g)}  \\
G_{\bf F}(n,N) & \stackrel {\quad \varphi \quad }
{\longrightarrow} & M_{\bf F}(N).  \end{array} $$

Let $\triangle =(d+\delta)^2$ be the Laplacian on $G_{\bf F}(n,N)$
with respect to the metric $2\Phi$. For any $A\in M_{\bf F}(N)$,
$f(\pi)=\langle \varphi(\pi), A\rangle $ is a function on the
Grassmann manifold $G_{\bf F}(n,N)$. As is well-known, we have
$$\triangle  f=-cn(N-n)\langle H, A\rangle.$$

As show above, $M_{\bf F}(n,N)$ is in the hyperplane $\{ B\in
M_{\bf F}(N) \ | \ \langle B, I_N\rangle=n\}$ of Euclidean space
$M_{\bf F}(N)$. For any vector $A$ parallel to this hyperplane, we
have
$$\langle A, I_N\rangle= \langle A,  \sum e_i\bar e_i^t\rangle +
\langle A,  \sum e_\alpha\bar e_\alpha^t\rangle =0.$$ Then for
such $A$, we have
$$ \triangle  f  = cN\langle \sum\limits_i  e_i \bar e_i^t ,
A\rangle= cNf.$$ We have proved the following

{\bf Theorem 2.4} \ The imbedding $\varphi: \ G_{\bf F}(n,N) \to
\{ B\in M_{\bf F}(N) \ | \  \mbox{tr} B=n\}$ is formed by the
eigenfunctions of Laplacian $\triangle$ on $G_{\bf F}(n,N)$ with
eigenvalue $cN$.

For any $A= \sum e_i\bar e_i^t\in M_{\bf F}(n,N)$, let $\widetilde
A= I_N- 2A=I_N- 2\sum e_i\bar e_i^t$. It is easy to see that
$\overline{\widetilde A}^t=\widetilde A, \ \widetilde A^2=I_N.$
Then $A\to \widetilde A$ gives a map $\psi: \ M_{\bf F}(n,N)\to
G_{\bf F}(N).$

It is interesting  to note that when ${\bf F}={\bf R}$ be the real
numbers, the imbedding of $G_{\bf R}(n,N)$ in orthogonal group
$O(N)$ defined above can be obtained by using Clifford algebra.
Let $C\ell_N$ be the Clifford algebra associated to the Euclidean
space ${\bf R}^N$ and $Pin(N)$ be the Pin group. Any unit vector
$v$ of ${\bf R}^N$ defines a reflection $f_v$ on ${\bf R}^N$:
$$f_v(e)= v\cdot e\cdot v=e - 2(e,v)v, \ \ \forall
e\in {\bf R}^N,$$ where `$\cdot$' denotes the Clifford product.
With the standard basis of ${\bf R}^N$, the map $f_v$ can be
represented by matrix $I_N-2v v^t\in O(N)$.

Let $\widetilde G_{\bf R}(n,N)$ be the oriented Grassmann
manifold. For any $\widetilde \pi\in \widetilde G_{\bf R}(n,N)$,
we choose an oriented orthonormal basis $e_1,\cdots, e_n$ of
$\widetilde \pi$. Note that $f_{e_i}f_{e_j}=f_{e_j}f_{e_i}$ for
any $i,j$. Then $e_1\cdot e_2\cdots e_n\in Pin(N)$ and by the maps
$\widetilde G_{\bf R}(n,N)\to Pin(N) \stackrel {Ad}\longrightarrow
O(N)$, we have a map
$$\widetilde G_{\bf R}(n,N)\to  O(N), \ \ \ \widetilde \pi \to
f_{e_1}f_{e_2}\cdots f_{e_n}, $$
\begin{eqnarray*}
& f_{e_1}f_{e_2}\cdots f_{e_n} & = (I_N-2e_1
e_1^t)(I_N-2e_2 e_2^t)\cdots (I_N-2e_n e_n^t) \\
&& =I_N -2\sum\limits_{i=1}^n e_ie_i^t. \end{eqnarray*}

The map $\widetilde G_{\bf R}(n,N)\to  O(N)$ is an imbedding. As
$\det(I_N -2\sum\limits_{i=1}^n e_ie_i^t)=(-1)^n$ is constant, we
can imbed  real Grassmann manifold $G_{\bf R}(n,N)$ in $SO(N)$.

\vskip 1cm \centerline{\bf \S 3. The Morse functions on the
Grassmann manifolds}

\vskip 0.3cm

In [12], we have constructed many (degenerate or non-degenerate)
Morse functions  on the real oriented Grassmann manifolds by using
calibrations. In the following we construct Morse functions on the
Grassmann manifolds $G_{\bf F}(n,N)$.

For any $A\in M_{\bf F}(N)$, \ the map $f: \ G_{\bf F}(n,N) \to
{\bf R}, \ f(\pi) = \ \langle \varphi(\pi), A\rangle, $ defines a
function on Grassmann manifolds $G_{\bf F}(n,N)$. $f$ is a Morse
function for almost every vector $A\in M_{\bf F}(N)$. But in
general, it is difficult to find such element. First, we give some
known results.

Let $E_{AB}$ be the elements in $M_{\bf F}(N)$, $A\geq B$, where
the entries in row $A$, \ column $B$ and  row $B$, column $A$ are
$1$, the others are zero.

When ${\bf F}={\bf C}, \ n=1, \ G_{\bf C}(1,N)\approx {\bf
C}P^{N-1}$ is the complex projective space. Let $A=\sum\limits_{A}
c_AE_{AA}\in M_{\bf C}(N), \ c_1>c_2>\cdots
>c_N>0$. For any $\pi\in G_{\bf C}(1,N), \ \varphi(\pi)=e_1\bar e_1^t, \
e_1=(z_1,z_2,\cdots,z_N)^t, \ \sum\limits_{A} |z_A|^2=1$, then
$$f(\pi)=\langle \varphi(\pi), A\rangle =\sum\limits_{A} c_A|z_A|^2.$$
As is well-known ([9]), $f$ is a perfect Morse function on ${\bf
C}P^{N-1}$.

Similar results hold for the real projective space $G_{\bf
R}(1,N)\approx {\bf R}P^{N-1}$ and the quaternion  projective
space $G_{\bf H}(1,N)\approx {\bf H}P^{N-1}$. In the real case,
the functions are not perfect.

The map $f: \ G_{\bf F}(n,N) \to {\bf R}$, $f(\pi)=\langle
\varphi(\pi), E_{11}\rangle, $ defines a function on $G_{\bf
F}(n,N)$. To study this function  we define two submanifolds of
$G_{\bf F}(n,N)$. Let $G_{\bf F}(n-1,N-1)$ be a submanifold of
$G_{\bf F}(n,N)$ such that every element of $G_{\bf F}(n-1,N-1)$
contains the vector $\tilde e_1 = (1,0,\cdots,0)^t\in {\bf
F}^{N}$. Let ${\bf F}^{N-1}=\{u=(0,u_2,\cdots,u_N)^t\in {\bf
F}^N\}$ be a subspace of ${\bf F}^{N}$ and  $G_{\bf F}(n,N-1)$ be
a submanifold of $G_{\bf F}(n,N)$ generated by the $n$-dimensional
subspaces of ${\bf F}^{N-1}$.

{\bf Theorem 3.1} \ The function  $f: \ G_{\bf F}(n,N) \to {\bf
R}$ is a degenerate Morse function on $G_{\bf F}(n,N)$, where
$f(\pi)=\langle \varphi(\pi), E_{11}\rangle $. The critical
submanifolds are $f^{-1}(0)=G_{\bf F}(n,N-1)$ and
$f^{-1}(1)=G_{\bf F}(n-1,N-1)$ with indices $0$ and $c(N-n)$
respectively.

{\bf Proof}. \ Let $e_1,\cdots,e_n, e_{n+1},\cdots,e_N$ be
orthonormal frame fields
 on ${\bf F}^N$ such that $G_{\bf F}(n,N)$ be generated by $e_1,\cdots,e_n$ locally.
We have $f(\pi)=\sum\limits_{i=1}^n x_{i1}\bar x_{i1},$ where
$e_i=(x_{i1},\cdots,x_{iN})^t$. Then, $0\leq f \leq 1$ and $\pi$
is a critical point of function $f$ if and only if
$$df=\sum\limits_{i,\alpha} \ \langle  e_\alpha\omega_i^\alpha \bar e_i^t +  e_i\bar
\omega_i^\alpha \bar e_\alpha^t, E_{11}\rangle =0.$$

We prove the theorem for the real case, the other cases can be
proved similarly, see the proof of Theorem 3.5. By Lemma 2.2,
$df=0$ if and only if
$$\langle  e_\alpha e_i^t +  e_i e_\alpha^t, E_{11}\rangle =0$$
for any $i,\alpha$. For $e_A=(x_{A1},\cdots,x_{AN}), \ A=1,\cdots,
N,$ we can assume  $x_{i1}=0$ for $i>1$ and $x_{\alpha 1}=0$ for
$\alpha
> n+1$. Obviously, we have $x_{11}^2 + x_{n+1 \ 1}^2=1$. Then
$$\langle  e_\alpha e_i^t +  e_i e_\alpha^t,
E_{11}\rangle =2x_{11}x_{n+1 \ 1}\delta_{1i}\delta_{\alpha n+1},$$
and the point $\pi\in G_{\bf R}(n,N)$ is a critical point if and
only if $x_{11}=0$ or $x_{11}=1$.

Let ${\bf R}^{N-1}=\{u=(0,u_2,\cdots,u_N)^t\in {\bf R}^N\}$ be a
subspace of ${\bf R}^{N}$ and  $G_{\bf R}(n,N-1)$ be a submanifold
of $G_{\bf R}(n,N)$ generated by the $n$-dimensional subspace of
${\bf R}^{N-1}$. Let $G_{\bf R}(n-1,N-1)$ be a submanifold of
$G_{\bf R}(n,N)$ such that every element of $G_{\bf R}(n-1,N-1)$
contains the vector $\tilde e_1 = (1,0,\cdots,0)^t\in {\bf
R}^{N}$. It is easy to see that $f^{-1}(0)=G_{\bf R}(n,N-1)$ and
$f^{-1}(1)=G_{\bf R}(n-1,N-1)$.

Now we show that the critical submanifolds $f^{-1}(0)$ and
$f^{-1}(1)$ of $f$ are non-degenerate and compute their indices.

On $f^{-1}(0)=G_{\bf R}(n,N-1), \ x_{i1}=0, \ i=1,\cdots,n, \
\tilde e_{n+1}=(1,0,\cdots,0)^t,$ then the tangent space of
$f^{-1}(0)$ is generated by
$$e_ie_\alpha^t+ e_\alpha e_i^t, \ \ \alpha\not= n+1.$$
On $f^{-1}(1)=G_{\bf R}(n-1,N-1), \ \tilde e_1 =(1,0,\cdots,0)$, \
$G_{\bf R}(n-1,N-1)$ is generated by $\tilde e_1,e_2,\cdots,e_n$,
then the tangent space of $f^{-1}(1)$ is generated by
$$e_ie_\alpha^t+ e_\alpha e_i^t, \ \ i\not= 1.$$

By  simple computation, on the critical submanifolds, we have
$$d^2f=-\sum \omega_j^\alpha\omega_i^\alpha\langle e_j  e_i^t +
 e_i e_j^t, E_{11}\rangle  + \sum \omega_i^\alpha\omega_i^\beta \langle
 e_\alpha e_\beta^t +
 e_\beta e_\alpha^t, E_{11}\rangle .$$
Then
$$d^2f|_{f^{-1}(0)}= 2\sum \omega_i^{n+1}\omega_i^{n+1} \langle
 e_{n+1} e_{n+1}^t, E_{11}\rangle  = 2\sum \omega_i^{n+1}\omega_i^{n+1},$$
$$d^2f|_{f^{-1}(1)}=-2\sum \omega_1^\alpha\omega_1^\alpha\langle e_1  e_1^t, E_{11}\rangle  =
-2\sum \omega_1^\alpha\omega_1^\alpha.$$ By Lemma 2.2, the
critical submanifolds of  $f$ are all non-degenerate. These
complete the proof of the theorem. \ \ \ $\Box$

By Morse theory, it can be shown that every differentiable
manifold has the homotopy type of a CW complex.  As in [4] or [5],
let
$$0\leq a_0\leq a_1\leq \cdots \leq a_{n}\leq N-n$$ be a sequence of
integers. These give a CW complex structure on the Grassmann
manifold $G_{\bf F}(n,N)$. For every such $(a_0, a_1, \cdots,
a_{n})$, there is one cell of dimension $c(a_0+ a_1+ \cdots +
a_{n})$. The homologies of the Grassmann manifold can be computed
by means of the Schubert varieties (see [4] or [5]). There is a
close relation between Theorem 3.1 and   the Schubert varieties:

The elements $(a_0, a_1, \cdots, a_{n})$ can be divided into two
classes:

(1) $(a_0, a_1, \cdots, a_{n})$, where $a_n \leq N-n-1$;

(2) $(a_0, a_1, \cdots, a_{n-1}, N-n)$, where $a_{n-1} \leq N-n$.

Let $\{\rho_1,\rho_2,\rho_3\}$ be a partition of unity on $[0,1]$
such that $\mbox {supp} (\rho_1)\subset [0,\frac 14],\\ \mbox
{supp} (\rho_2)\subset [\frac 18,\frac 78], \ \mbox {supp}
(\rho_3)\subset [\frac 34,1]$ and $\frac {d\rho_1}{dt} \leq 0, \
\frac {d\rho_3}{dt} \geq 0$. Let $h_1\leq 0,h_2\geq 0$ be two
non-degenerate Morse functions on $f^{-1}(0)$ and $f^{-1}(1)$
respectively. The functions $h_1, h_2$ can be viewed as functions
on neighborhoods of $f^{-1}(0), f^{-1}(1)$ in $G_{\bf F}(n,N)$
respectively, they are constants on the trajectories of grad$(f)$.
Define a function
$$\tilde f = \tilde \rho_1 (h_1+ f) +\tilde \rho_2 f +  \tilde \rho_3(h_2+
f)= \tilde \rho_1 h_1 + \tilde \rho_3h_2+ f,$$ where $\tilde
\rho_i=\rho_i\circ  f, \ i=1,2,3.$

 {\bf Theorem 3.2} \ $\tilde f: \ G_{\bf F}(n,N) \to {\bf R}$ is a non-degenerate Morse
function  and the critical points are that of $h_1$ and $h_2$. If
$p$ is a critical point of $h_2$ with index $k$, then the index of
$p$ is $k+c(N-n)$ with respect to the function $\tilde f$; if $q$
is a critical point of $h_1$,  then the indices of $q$ with
respect to the functions $h_1$  and $\tilde f$ are the same.

For the proof,  see [12].

When ${\bf F}={\bf C}$ or ${\bf H}$, we can choice perfect Morse
functions $h_1,h_2$ on  $G_{\bf F}(n,N-1)$ and $G_{\bf
F}(n-1,N-1)$ respectively. Then $\tilde f$ is also a perfect Morse
function. Let $P_t(M)$ be the Poincar\'{e} polynomial for a
manifold $M$.

{\bf Corollary 3.3} \ For ${\bf F}={\bf C}$ or ${\bf H}$, the
Poincar\'{e} polynomial of $G_{\bf F}(n,N)$ can be represented by
$$P_t(G_{\bf F}(n,N))=P_t(G_{\bf F}(n,N-1)) +t^{c(N-n)}P_t(G_{\bf F}(n-1,N-1)).$$

For example, by  simple computation, we have
$$P_t(G_{\bf C}(1,N))=1+t^2+t^4+\cdots +t^{2(N-1)},$$
$$P_t(G_{\bf C}(2,5)) = 1+ t^2+2t^4+2t^6+2t^8+t^{10}+t^{12},$$
$$P_t(G_{\bf C}(2,7)) = 1+ t^2+2t^4+2t^6+3t^8+3t^{10}+3t^{12}
+2t^{14}+2t^{16} +t^{18}+t^{20},$$
$$P_t(G_{\bf C}(2,8)) = P_t(G_{\bf C}(2,7))+t^{12}P_t(G_{\bf C}(1,7)),$$
$$P_t(G_{\bf C}(3,7))=(1+t^6)P_t(G_{\bf C}(2,5))+ t^8P_t(G_{\bf C}(2,6)),$$
$$P_t(G_{\bf C}(5,10)) = (1+ t^{10})[t^{20}P_t(G_{\bf C}(2,7))+
(1+t^8+t^{10})P_t(G_{\bf C}(3,7))].$$

By $P_t(G_{\bf F}(n,N))=P_t(G_{\bf F}(N-n,N))$, we have
$$(t^{cn}-1)P_t(G_{\bf F}(n,N-1))=(t^{c(N-n)}-1)P_t(G_{\bf
F}(n-1,N-1)).$$

Now we study the trajectories of gradient vector field of the
degenerate Morse function $f: \ G_{\bf F}(n,N) \to {\bf R}$
defined in Theorem 3.1. The gradient of the function $f: \ G_{\bf
F}(n,N) \to {\bf R}$ is
$$\mbox {grad} (f) =\frac 12 \sum\limits_{\tau}
\langle \xi_\tau, E_{11}\rangle\xi_\tau,$$ where the tangent
vectors $\xi_\tau$ are defined as in Lemma 2.1. For any $\pi\in
G_{\bf F}(n,N) - (f^{-1}(0)\cup f^{-1}(1))$, let
$e_i=(x_{i1},x_{i2},\cdots,x_{iN})^t, \ i=1,\cdots,n,$ be an
orthonormal basis of $\pi$ such that $0<x_{11}< 1, \ x_{i1}=0$ for
$i>1$. Denote $e_1(t)= (\cos t, (\sin t )x)^t$, where
$x=(x_{12},\cdots,x_{1N})/\sqrt{|x_{12}|^2+\cdots + |x_{1N}|^2}$.
Then there is $t_0\in (0,\frac {\pi}2)$ such that $e_1(t_0)=e_1$.
Let $\gamma(t)$ be a curve in $G_{\bf F}(n,N)$ generated by
orthonormal vectors $e_1(t),e_2,\cdots, e_n$. Then
$$f(\gamma(t))=\cos^2t, \ \ \ \frac {df(\gamma(t))}{dt} =-\sin 2t,$$
$$\gamma(0)\in f^{-1}(1)=G_{\bf F}(n-1,N-1), \ \gamma(\frac {\pi}2)\in
f^{-1}(0)=G_{\bf F}(n,N-1).$$ Note that $\dim \gamma(0)\cap
\gamma(\frac {\pi}2)=c(n-1)$. Along the curve $\gamma(t)$, let
$$e_{n+1}(t) =(-\sin t, (\cos t ) x)^t, \ \
e_\alpha =(0,x_{\alpha 2},\cdots,x_{\alpha N})^t, \ \ \alpha
>n+1,$$
be orthonormal complement of the vectors $e_1(t),e_2,\cdots, e_n$
in ${\bf F}^N$. Therefore
$$\mbox {grad} (f)|_{\gamma} = -\frac 12\sin 2t
(e_{n+1}(t) \bar e_1^t(t) + e_1(t)\bar e_{n+1}^t(t))= -\frac
12\sin 2t\frac {d\gamma}{dt}.$$ This shows that the curve $\gamma$
is a trajectory of the vector field $\mbox {grad} (f)$ on $G_{\bf
F}(n,N)$.

It is also easy to see that the vector  $\frac {d\gamma}{dt}(0)$
is normal to $f^{-1}(1)$ and the vector  $\frac
{d\gamma}{dt}(\frac {\pi}2)$ is normal to $f^{-1}(0)$. Let ${\bf
F}P^{N-n-1}=G_{\bf F}(1, N-n)$ be a subspace of $f^{-1}(0)$ such
that $e_2,\cdots,e_n\in \pi$ for any $\pi\in {\bf F}P^{N-n-1}$.
Let ${\bf F}P^{n-1}=G_{\bf F}(1, n-1)$ be a subspace of
$f^{-1}(1)$, any $\pi\in {\bf F}P^{n-1}$ be generated by
$e_1=(1,0,\cdots,0)^t, \tilde e_2,\cdots,\tilde e_n$, where
$\tilde e_2,\cdots,\tilde e_n\in \gamma(\frac {\pi}2)$.

{\bf Theorem 3.4} \ The trajectories of $\mbox {grad} (f)$ give
the maps from ${\bf F}P^{N-n-1}$ to $\gamma(0)$ and ${\bf
F}P^{n-1}$ to $\gamma(\frac {\pi}2)$ respectively.

When $n=1$, these gives the following canonical cell decomposition
of the projective space ${\bf F}P^{N-1}$
$${\bf F}P^{0}\subset {\bf F}P^{1}\subset\cdots \subset  {\bf F}P^{N-2}\subset
{\bf F}P^{N-1}.$$

In the following we assume $n, N-n\geq 2$.

{\bf Theorem 3.5} \ Let $g: \ G_{\bf F}(n,N) \to {\bf R}, \
g(\pi)=\langle \varphi(\pi),E_{12}\rangle $. The function $g$ is a
degenerate Morse function with critical submanifolds
$g^{-1}(0)=G_{\bf F}(n-2, N-2)\bigcup G_{\bf F}(n, N-2), \
g^{-1}(-1)=\widetilde G_{\bf F}(n-1, N-2), \ g^{-1}(1)= G_{\bf
F}(n-1, N-2)$. The indices on $G_{\bf F}(n-2, N-2), G_{\bf F}(n,
N-2), \ \widetilde G_{\bf F}(n-1, N-2), \  G_{\bf F}(n-1, N-2)$
are $ c(N-n), cn, 0, c(N-1)$ respectively.

{\bf Proof}. \ Let $e_i=(x_{i1},x_{i2},\cdots,x_{iN})^t$, then,
$$g(\pi) = \mbox{Re} \ \sum\limits_i (x_{i1}\bar x_{i2} + x_{i2}\bar
x_{i1})=2\mbox{Re} \ \sum\limits_i x_{i1}\bar x_{i2}.$$ We can
assume $x_{i1}=0$ for $i>1$ and $x_{11}$  a real number, this
shows $-1\leq g(\pi)\leq 1$. The critical points of function $g$
are determined by
$$dg=\sum\limits_{i,\alpha} \ \langle  e_\alpha\omega_i^\alpha \bar e_i^t +  e_i\bar
\omega_i^\alpha \bar e_\alpha^t, E_{12}\rangle =0,$$ where
$\varphi(\pi)=\sum e_i\bar e_i^t$.

We prove the theorem for the  case of ${\bf F}={\bf H}$. $dg=0$ if
and only if
$$\langle  e_\alpha \bar e_i^t +  e_i \bar e_\alpha^t, E_{12}\rangle=0, \  \
\langle  e_\alpha i\bar e_i^t -  e_i i\bar e_\alpha^t,
E_{12}\rangle =0,$$
$$\langle  e_\alpha j\bar e_i^t -  e_ij \bar e_\alpha^t, E_{12}\rangle =0, \  \
\langle  e_\alpha k\bar e_i^t -  e_i k\bar e_\alpha^t,
E_{12}\rangle =0,$$ for any $i,\alpha$.  Obviously, we have
\begin{eqnarray*} & & \langle  e_\alpha
\bar e_i^t +  e_i \bar e_\alpha^t, E_{12}\rangle  \\
& & = \mbox {Re} \  (x_{\alpha 1}\bar x_{i2} + x_{\alpha 2}\bar
x_{i1} +x_{i1}\bar x_{\alpha 2}
+ x_{i 2}\bar x_{\alpha 1}) \\
& & =2\mbox {Re} \  [x_{\alpha 1}\bar x_{i2} + x_{\alpha 2}\bar
x_{i1}],
\end{eqnarray*}
and
\begin{eqnarray*} & & \langle  e_\alpha i\bar e_i^t -  e_i i\bar e_\alpha^t, E_{12}\rangle  \\
& & = \mbox {Re} \  [(x_{\alpha 1}i\bar x_{i2} -x_{i 2}i\bar
x_{\alpha 1}) +
( x_{\alpha 2}i\bar x_{i1}- x_{i1}i\bar x_{\alpha 2})] \\
& & =2\mbox {Re} \ [x_{\alpha 1}i\bar x_{i2}  +   x_{\alpha
2}i\bar x_{i1}]. \end{eqnarray*} Similarly,
$$\langle  e_\alpha j\bar e_i^t -  e_i j\bar e_\alpha^t, E_{12}\rangle
=2\mbox {Re} \ [x_{\alpha 1}j\bar x_{i2}  +   x_{\alpha 2}j\bar
x_{i1}],$$ $$\langle  e_\alpha k\bar e_i^t -  e_i k\bar
e_\alpha^t, E_{12}\rangle  =2\mbox {Re} \ [x_{\alpha 1}k\bar
x_{i2} + x_{\alpha 2}k\bar x_{i1}].$$ {\bf F}or any $u,v\in {\bf
H}, a\in \mbox{Im}{\bf H}$, the following holds
$$\mbox {Re} \  (ua\bar v)=\mbox {Re} \  (-u\overline {va})
=\mbox {Re} \  (-\bar uva).$$ These show $\pi$ is a critical point
of $g$ if and only if
$$x_{\alpha 1}\bar x_{i2} + x_{\alpha 2}\bar x_{i_1}=0, \ \ \mbox{for all} \ i,\alpha.$$

On the critical submanifolds, we have
$$d^2g=-\sum\limits_{i,j,\alpha} \langle e_j \bar \omega_j^\alpha\omega_i^\alpha \bar e_i^t +
 e_i \bar\omega_i^\alpha\omega_j^\alpha \bar e_j^t, E_{12}\rangle  + \sum\limits_{i,j,\alpha}
   \langle
 e_\alpha \omega_i^\alpha\bar \omega_i^\beta \bar e_\beta^t +
 e_\beta \omega_i^\beta\bar\omega_i^\alpha\bar e_\alpha^t, E_{12}\rangle .$$

(1) Let ${\bf H}^{N-2}=\{ (0,0,x_3,\cdots, x_N)^t\in {\bf H}^N\}$
be a subspace of ${\bf H}^N$ and $G_{\bf H}(n, N-2)=\{ \pi\in
G_{\bf H}(n,N) \ | \ \pi\subset {\bf H}^{N-2}\}$ be submanifold of
$G_{\bf H}(n,N)$. Then $G_{\bf H}(n, N-2)$ is a critical
submanifold of function $g$ and $ g|_{G_{\bf H}(n, N-2)} \equiv
0$. In this case, $x_{i1}=x_{i2}=0$ for $i=1,\cdots,n$, we can
assume
$$e_{n+1}=(1,0,0,\cdots, 0)^t, \
e_{n+2}=(0,1,0,\cdots, 0)^t.$$ Then
\begin{eqnarray*} & & d^2g|_{G_{\bf H}(n, N-2)}\\
& &  =2\mbox{Re} \ \sum\limits_{i} (\omega_i^{n+1}\bar
\omega_i^{n+2}+\omega_i^{n+2}\bar \omega_i^{n+1}) \\
& & =2\sum\limits_{i} [ (a_i^{n+1}+a_i^{n+2})^2 -
(a_i^{n+1}-a_i^{n+2})^2 +
(b_i^{n+1}+b_i^{n+2})^2 - (b_i^{n+1}-b_i^{n+2})^2 \\
& & \quad + (c_i^{n+1}+c_i^{n+2})^2 - (c_i^{n+1}-c_i^{n+2})^2 +
(d_i^{n+1}+d_i^{n+2})^2 - (d_i^{n+1}-d_i^{n+2})^2],
\end{eqnarray*} where $\omega_i^\alpha =a_i^\alpha+ib_i^\alpha
+jc_i^\alpha+kd_i^\alpha.$ As in the proof of Theorem 3.1, we can
show that the critical submanifold $G_{\bf H}(n, N-2)$ is
non-degenerate with index  $4n$.

(2) Let $\tilde e_1=(1,0,\cdots,0)^t, \ \tilde
e_2=(0,1,0,\cdots,0)^t\in {\bf H}^N$ and $G_{\bf H}(n-2, N-2)=\{
\pi\in G_{\bf H}(n,N) \ | \ \tilde e_1,\tilde e_2\in \pi\}$ be a
submanifold of $G_{\bf H}(n,N)$. Then we have $x_{\alpha
1}=x_{\alpha 2}=0, \ \alpha=n+1,\cdots,N,$ for any $\pi\in G_{\bf
H}(n-2, N-2)$. Therefore $g|_{G_{\bf H}(n-2, N-2)} \equiv 0$ and
$G_{\bf H}(n-2, N-2)$ is a critical submanifold of $g$,
$$d^2g|_{G_{\bf H}(n-2, N-2)}
=-2\mbox{Re} \ (\sum\limits_{\alpha} \bar
\omega_1^\alpha\omega_2^\alpha  +
\bar\omega_2^\alpha\omega_1^\alpha ).$$ The critical  submanifold
$G_{\bf H}(n-2, N-2)$ is non-degenerate with index $4(N-n)$.

(3) Now we study  the case of the numbers $x_{i1},x_{i2}$ are not
all zeros and so are  the numbers $x_{\alpha 1},x_{\alpha 2}$.
Assuming $x_{i1}=0$ for $i>1$, \ $x_{j2}=0$ for $j>2$; \
$x_{\alpha 1}=0$ for $\alpha> n+1$, \ $x_{\beta 2}=0$ for
$\beta>n+2$. Then the conditions $x_{\alpha 1}\bar x_{i2} +
x_{\alpha 2}\bar x_{i_1}=0$ become
$$\frac {\bar x_{11}}{\bar x_{12}} =\frac {0}{\bar x_{22}}=-\frac {x_{n+1 \ 1}}{x_{n+1 \ 2}}
=-\frac {0}{x_{n+2 \ 2}}.$$ If $x_{11}=0$, we can assume
$x_{22}=0$. Then if $x_{11}=0$, we  have $x_{12}=0$. Therefore
$x_{11}\not =0,x_{12}\not=0$ in this case. Similarly, $x_{n+1 \
1}\not =0,x_{n+1 \ 2}\not=0$. Hence $x_{22}=x_{n+2 \ 2}=0$. By
$$\frac {\bar x_{11}}{\bar x_{12}} =-\frac {x_{n+1 \ 1}}{x_{n+1 \ 2}},
 \ \ |x_{11}|^2+|x_{n+1 \ 1}|^2=1, \ |x_{12}|^2+|x_{n+1 \ 2}|^2=1,$$
and the vectors $e_1\perp e_{n+1}$,
 we have
$$e_1=(\frac { \ 1}{\sqrt 2},\pm \frac { \ 1}{\sqrt 2},0,\cdots,0)^t,
\ \ e_i=(0,0,x_{i3},\cdots,x_{iN})^t,\ \ i>1.$$

Let $G_{\bf H}(n-1,N-2)$ be the subset of  $\pi\in G_{\bf H}(n,N)$
which is generated by $e_1=(\frac { \ 1}{\sqrt 2},\frac { \
1}{\sqrt 2},0,\cdots,0)^t, \ \ e_i=(0,0,x_{i3},\cdots,x_{iN})^t,\
\ i>1$. Similarly, Let $\widetilde G_{\bf H}(n-1,N-2)\subset
G_{\bf H}(n,N)$ be the subset of $\pi$ which is generated by
$\tilde e_1=(\frac { \ 1}{\sqrt 2}, - \frac { \ 1}{\sqrt
2},0,\cdots,0)^t, \ \ e_i=(0,0,x_{i3},\cdots,x_{iN})^t,\ \ i>1$.
By construction, $G_{\bf H}(n-1,N-2)$ and $\widetilde G_{\bf
H}(n-1,N-2)$ are critical submanifolds of function $g$, \
$$g|_{G_{\bf H}(n-1,N-2)}\equiv 1, \ \ g|_{\widetilde G_{\bf H}(n-1,N-2)}\equiv -1.$$
By our assumption, $g(\pi) =2\mbox{Re} \ (x_{11}\bar x_{12})$ and
$|x_{11}|^2+ |x_{12}|^2\leq 1$, this shows
$$g^{-1}(1)=G_{\bf H}(n-1,N-2), \ \ g^{-1}(-1)=\widetilde G_{\bf H}(n-1,N-2).$$

On $G_{\bf H}(n-1,N-2)$, we can set  $e_{n+1}=(\frac { \ 1}{\sqrt
2}, -\frac { \ 1}{\sqrt 2},0,\cdots,0)^t$; on $\widetilde G_{\bf
H}(n-1,N-2)$, we can set   $\tilde e_{n+1}=(\frac { \ 1}{\sqrt 2},
\frac { \ 1}{\sqrt 2},0,\cdots,0)^t$. Then we have
$$d^2g|_{g^{-1}(-1)}=\mbox{Re} \  [\sum\limits_{\alpha}
 \bar\omega_1^\alpha\omega_1^\alpha  +
\sum\limits_{i} \omega_i^{n+1}\bar\omega_i^{n+1} ],$$
$$d^2g|_{g^{-1}(1)}=-\mbox{Re} \  [\sum\limits_{\alpha}
 \bar\omega_1^\alpha\omega_1^\alpha +
\sum\limits_{i} \omega_i^{n+1}\bar\omega_i^{n+1} ].$$

As in the proof of Theorem 3.1, we can show that the critical
submanifolds $g^{-1}(-1)$ and $g^{-1}(1)$ are non-degenerate with
indices $0$ and $4(N-1)$ respectively. \ \ \ $\Box$

As Corollary 3.3, we have

{\bf Corollary 3.6} \ For ${\bf F}={\bf C}$ or ${\bf H}$, the
Poincar\'{e} polynomial of $G_{\bf F}(n,N)$ can be represented by
\begin{eqnarray*} & P_t(G_{\bf F}(n,N))= &
t^{cn}P_t(G_{\bf F}(n,N-2)) +t^{c(N-n)}P_t(G_{\bf F}(n-2,N-2)) \\
& & + (1+ t^{c(N-1)})P_t(G_{\bf F}(n-1,N-2)),
\end{eqnarray*} where $n, N-n\geq 2$.

\vskip 1cm \centerline{\large \bf References}

\vskip 0.3cm {\small

 \noindent [1]  \  Bott, R., \ Nondegenerate
critical manifolds, Ann. of Math., {\bf 60}(1954), 248-261.

\noindent [2] \   Bott, R., \ Lectures on Morse theory, Bull.
Math. Soc. {\bf 7}(1982), 331-358.

\noindent [3] \ Chen, W. H.,    \ The differential geometry of
Grassmann manifold as submanifold, Acta Math., Sinica, {\bf
31A}(1988), 46-53.

\noindent [4]  \  Chern, S. S.,  \    Complex manifolds without
potential theory,  Springer-Verlag, New York, 1979.

\noindent [5] \ Chern, S. S., \ Topics in differential geometry,
Inst. for Adv. Study, Princeton, 1951.

\noindent [6] \ Chern, S. S.,  do Carmo M. and Kabayashi, S.,  \
Minimal submanifolds of a sphere with second fundamental form of
constant length, Functional Analysis and Related Topics,
Springer-Verlag(1970), 59-75.

\noindent [7] \  Helgason, S.,  \ Differential geometry, Lie
groups, and symmetric spaces, Academic Press, New York, 1978.

\noindent [8]  \ Kabayashi, S. and Nomizu, K., \ Foundations of
differential geometry, vol. 2,  Interscience Publishers, New York,
1969, \ 159-161.

\noindent [9] \ Milnor, J., \ Morse theory, Princeton University
Press, 1963.

\noindent [10] \ Takahashi, T., \ Minimal immersions of Riemannian
manifolds,  J. Math. Soc, Japan, {\bf 18}(1966), 380-385.

\noindent [11]  Takeuchi, M and Kobayashi, S.,  Minimal imbeddings
of R-spaces, J. Differential Geometry,  2 (1968), 203-215.

\noindent [12] \ Zhou, J. W., \ Morse Functions on Grassmann
Manifolds, Proc. of the Royal Soc. of Edinburgh, 135A(2005),
209-221.

\end{document}